\documentclass{article}

\usepackage{arxiv}

\usepackage[utf8]{inputenc} 
\usepackage[T1]{fontenc}    
\usepackage{hyperref}       
\usepackage{url}            
\usepackage{booktabs}       
\usepackage{amsfonts}       
\usepackage{nicefrac}       
\usepackage{microtype}      
\usepackage{lipsum}
\usepackage{graphicx}
\usepackage{amsmath,amsfonts}
\usepackage{float}

\graphicspath{ {./images/} }

\providecommand{\edot}[1]{\ensuremath{\cdot10^{#1}}} 

\title{Evaluation of a Fractional-Calculus-based Numerical Approach to solve Ordinary Differential Equations}

\author{
 Sergio F. Yapur \\
  Facultad de Ingenier\'ia Qu\'imica \\
  Universidad Nacional del Litoral\\
  Santiago del Estero 2829 (3000) Santa Fe \\
  \texttt{syapur@fiq.unl.edu.ar} \\
}

\begin{document}
\maketitle
\begin{abstract}
This article examines a new approach to solving ordinary differential equations based on Fractional-Calculus theory. Poisson and Sturm-Liouville-type problems are studied, together with different boundary conditions. Each case is analyzed and compared concerning the Finite-Difference method outcome.

\end{abstract}


\section{Introduction}

Several numerical methods for solving ordinary differential equations use numerical differentiation techniques. These approximations generally consist of discrete rational functions. As such, the quality of the approximation relies on the discrete increment value $h>0$. An unsuitable magnitude of h could lead to numerical instabilities since the small values of $h$ needed to reduce truncation error also cause the round-off error to grow \cite{Burden}.

On the other hand, the Fractional Differential Equations (FDE) theory continues to consolidate through an expanding collection of books and publications dealing with resolution techniques. However, to the best of my knowledge, there are no attempts to improve conventional Ordinary Differential Equations (ODE) theory starting from the emerging FDE. The new theoretical framework enables us to revisit classical problems from the numerical perspective, at least in principle.

This article explores an application of FDE by building a numerical scheme that tempers the strong dependency with respect to the discrete increment h. In general terms, the idea is to add a new dimension over which discretization is possible. A non-integer order of integration will define this dimension.

As mentioned above, there is a sensitivity problem regarding the increment $h$ when using numerical differentiation. To lower the approximation error, $h$ should be smaller whenever the slope of the associated function is greater. In fact, this is the working principle of adaptive methods. By adjusting $h$ locally, it is possible to minimize the approximation error for a given computational cost. Alternatively, would it be possible to keep a relatively large value of $h$ by adjusting the rate of change? This workaround is feasible by taking a derivative of fractional order $\epsilon$, with $0<\epsilon<1$. The main idea is to exploit fractional derivatives to maximize the increment value h while keeping a given error. This article explores the computational aspects of this approach.

It is desirable to have a reference method to compare scheme performances. The Finite Difference Method (FDM) is adopted in this work, primarily because of its low computational cost. This feature makes the FDM a preferential method to solve ODE on simple domains.

\subsection{Brief Introduction to Fractional Calculus}

The so-called Fractional Calculus deals with operators that express varying degrees of integral or derivative behavior, depending on the value of the order. Because of this polyvalence, they are called differintegral operators. The order $\alpha$ is complex, being equal to the usual differentiation whenever $\alpha \in \mathbb{R}^+$, as well as to the usual integration whenever $\alpha \in \mathbb{R}^-$.

It is considered that the theory originated in a letter from l'H\^opital to Leibniz in 1695. He wondered about the implications of having a derivative $d^n f/dx^n$ of a suitable function f with $n=1/2$. By that time, such a theoretical artifact was possibly too advanced, and this operator was unknown. But the question remained, and later developments over the following centuries were due to mathematicians such as Abel, Liouville, Letnikov, among many others. A detailed list of foundational works is available in the book of \cite{Spanier1974}. From 1970, this realm took new impetus with emerging applications in science and engineering. Some examples are anomalous diffusion, Brownian modeling, and control theory. The book of \cite{Baleanu} describes this fast expansion of applications. 

Even when there are a few definitions of differintegral operators, that of Riemann-Liouville (RL) is the most widely used. Let $f$ be a piecewise continuous function in $(0,\infty)$, integrable over any finite interval $[0,\infty)$. Given $c \in \mathbb{R}_0^{+}$ and $x>c$ the RL operator is:

\begin{equation}\label{Dfrac}
{}_{c} D_{x}^{\alpha} f(x) = 
\left\{ 	
	\begin{array}{ll}
	\displaystyle \frac{1}{\Gamma(-\alpha)}  \int_{c}^{x} \frac{f(\zeta) }{(x-\zeta)^{1+\alpha}}d\zeta, &  \mbox{si $\alpha \in \mathbb{R}^-$,} \\[12pt]
	f(x)  & \mbox{si $\alpha=0$,} \\[6pt]
	\displaystyle \frac{d^{\lceil\alpha\rceil}}{dx^{\lceil\alpha\rceil}} \, {}_{c}D_{x}^{\alpha-\lceil\alpha\rceil} f(t) & \mbox{si $\alpha \in \mathbb{R}^+$}
	\end{array}
\right.
\end{equation}

where $\lceil \alpha \rceil$ denotes the ceiling function of $\alpha$, and $\Gamma$ is the usual gamma function,

\begin{equation}
\Gamma(x)= \int_{0}^{\infty} e^{-y} y^{x-1} dy, \quad \mbox{si $x\in \mathbb{R}^+$.}
\end{equation}

Notice that definition \ref{Dfrac} is recursive when $\alpha>0$. Although this is not the most general definition, it can be easily extended to the case when  $\alpha \in \mathbb{C}$. In any case, the definition \ref{Dfrac} will be enough for this work. From now on $c=0$ is considered.

Regarding numerical applications, the Gr\"{u}nwald-Letnikov definition gives a useful expression for the same kind of operator.

\begin{equation}\label{DGL}
	{}_{c} D_{x}^{\alpha} f(x) = \lim\limits_{h \rightarrow 0} \frac{1}{h^\alpha} \sum_{j=0}^{[(x-c)/h]}(-1)^j \binom{\alpha}{j}f(x-jh),
\end{equation}

where $[.]$ denotes the integer part function, and the generalized binomial coefficient is used for non-integer arguments. Notice that the sum range is infinite. All these operators are, in fact, non-local whenever $\alpha \notin \mathbb{N}$. The infinite summation does not represent an issue from a numerical standpoint, as there are ways to approximate this operator. One of the most popular is the Principle of Finite Memory \cite{Baleanu}.

Other important definitions for the differintegral operator are due to Caputo and Hadamard. They have their strengths and are described in the bibliography \cite{Petras2011, Baleanu}. While the definitions \ref{Dfrac} and \ref{DGL} are equivalent \cite{Spanier1974}, the Caputo version cannot be interchanged with either the RL or the GL definitions in general. The main difference lies in the kernel structure. However, the three definitions are equivalent for differential problems with null initial conditions \cite{Rahimy}. Hereinafter, the RL definition will be used.

Some typical characteristics of ordinary derivatives are lost when using the fractional operator \ref{Dfrac}. For instance, ${}_{c} D_{x}^{\alpha} f(x) = 0$ does not imply, in general, that $f$ is a constant function when $\alpha \ne  1$. Fortunately, other properties do remain. A remarkable one is the semigroup property regarding the orders of similar differintegral operators:

\begin{equation}\label{semigrupo}
{}_{c} D_{x}^{\alpha} f(x) {}_{c} D_{x}^{\beta} f(x) = {}_{c} D_{x}^{\alpha+\beta} f(x) 
\end{equation}

This property, proved by Letnikov \cite{Spanier1974}, will turn out fundamental to building the technique proposed later on in this paper.

\section{CASE STUDIES}

This work presents four case studies. They all are classical EDO, whose structures are widely known for the variety of their applications. In each case, the objective is to analyze the numerical scheme performance.

The first case deals with the following one-dimensional ODE with Dirichlet boundary conditions.

\begin{equation}\label{caso1}
\mbox{Caso 1}\left\{ 
\begin{array}{rcl}
u''(x) &=& -20 e^{-10(x-0.7)^2}, \quad 0 \le x \le 1 \\
u(0) &=& a_1,  \\
u(1) &=& b_1	
\end{array}
\right.
\end{equation}

Where $a_1,b_1 \in \mathbb{R}$ are constants whose values will be stated later on. The second case corresponds to the same type of ODE with Dirichlet conditions at one end of the interval and Robin conditions at the other end.

\begin{equation}\label{caso2}
\mbox{Caso 2}
\left\{ 
\begin{array}{rcl}
u''(x) &=& -20 e^{-10(x-0.7)^2}, \quad 0 \le x \le 1 \\
u(0) &=& a_2,  \\
u'(1) + b_2 u(1) &=& c_2 	 	 
\end{array}
\right.
\end{equation}

Here, $a_2,b_2, c_2 \in \mathbb{R}$ are constants associated with case 2. The third case also deals with a Poisson-type equation with mixed boundary conditions, but the excitation term is substantially different from the previous cases.

\begin{equation}\label{caso3}
\mbox{Caso 3}
\left\{ 
\begin{array}{rcl}
u''(x) &=& 	-x(1-\sin^2(100x)), \quad 0 \le x \le 1 \\
u(0) &=& a_3,  \\
u'(1) + b_3 u(1) &=& c_3	 	 
\end{array}
\right.
\end{equation}

Again, $a_3,b_3,c_3 \in \mathbb{R} $ are generic constants. Finally, the last case is slightly more general, being a one-dimensional Sturm-Liouville problem with Dirichlet boundary conditions.

\begin{equation}\label{caso4}
\mbox{Caso 4}
\left\{ 
\begin{array}{rcl}
u''(x) &=& 	2x(5 - u), \quad 0 \le x \le 1 \\
u(0) &=& a_4,  \\
u(1) &=& b_4		 	 
\end{array}
\right.
\end{equation}

As before, $a_4,b_4 \in \mathbb{R}$ are constants. Notice that in each case the domain is the closed interval $[0,1]$.

\section{METHODS}

\subsection{Basic Principle}

The semigroup property enables an integration of the ODE in a recursive manner by using the expression:

\begin{eqnarray}\label{iterIOFI}
f_i &=& {}_{c} D_{x}^{\alpha_i} f_{i-1}, i=1,...,n \\
f_1 &=& Q,		\nonumber
\end{eqnarray}

The symbol $f_i$ denotes an iterated function and $Q$ represents the seed, as well the right-hand side of the equation, such as $u''(x) = Q(x)$. Also, $-2\le \alpha_i < 0, \forall i$. Recursion is invoked $n$ times until the condition $\sum_{1}^{n} \alpha_i = 2$ is satisfied, meaning that the procedure reached the usual second-order derivative. By convenience, this kind of approach will be called Iterative Fractional Order Integration (IFOI).

\subsection{Numerical Schemes}

Several techniques may approach differintegral operators for IFOI methods. The most straightforward calculation uses the GL definition given by Eq. \ref{DGL}. This expression entails a natural numerical approximation \cite{Feliu2010} by replacing the limit value of $h$ with a small $h>0$. Also, rearranging the summation coefficients recursively allows greater computational efficiency. The resulting algorithm is studied in Case 1. For Case 2, the so-called \textit{product rectangle formula} \cite{Baleanu} helps to approximate Eq. \ref{Dfrac}. Finally, for Cases 3 and 4, a variant of the predictor-corrector method for FDE, which is originally due to Adams-Bashforth-Moulton and was later modified by Deng (ABMD), is implemented.

On the other hand, the FDM for lineal problems is the usual \cite{Larsson2009}. It will be the method of reference for cases 1, 2, and 3. For case 4, a variant for the non-linear problem is used \cite{Burden}.

The reported errors below are absolute. They approach the supremum norm like follows:

\begin{equation}
e = d_{\infty}(\tilde{u},u) = \| \tilde{u} - u \|_\infty \approx \max_{0 \le i \le n} \left| \tilde{u}(x_i) - u(x_i) \right|,
\end{equation}

here, $u$ is the analytical solution to the problem, $\tilde{u}$ is the numerical approximation and $x_0<...<x_n$ is a partition of $[0,1]$, with $x_0=0$ and $x_n=1$.

Boundary conditions are treated with the \textit{shooting method} \cite{Isaacson, Burden}. It facilitates the solution of a boundary value problem (BVP) by solving two initial value problems (IVP). The final solution to the BVP is expressed as a linear combination:

\begin{equation}\label{CLshooting}
u = u_1 + c\,u_2,
\end{equation}

where $u_1,u_2$ are the individual solutions to each IVP, and the coefficient $c$ is obtained from the particular forms of the boundary conditions. More precisely, $u_2$ is the solution to the homogeneous problem, while $u_1$ is a particular solution to the general problem. 

The \textit{shooting} technique works well with either Robin or Dirichlet types boundary conditions. Also, due to its analytical nature, it does not affect the performance of numerical methods regarding approximation error.

\subsection{Tools}

The computing software is Matlab\textsuperscript{\tiny\textregistered}, version 2018a. The hardware is a PC with a processor Intel i7-3770 @ 3.4 GHz and 8 Gb RAM.

Every algorithm is developed \textit{ad hoc} for this work unless stated otherwise. Also, all codes were validated against several examples.

\section{RESULTS}
\label{sec:Results}

\subsection{Case 1}

The GL definition approaches the differintegral operator in this problem. An iterative algorithm was taken from the work of from \cite{Feliu2010}. This code is attractive due to its simplicity and speed in comparison to others.

The constants to specify problem \ref{caso1} are set to $a_1=-3$ and  $b_1=-2$.

\begin{figure}[ht]
	\centering
	\includegraphics[width=1.0\textwidth]{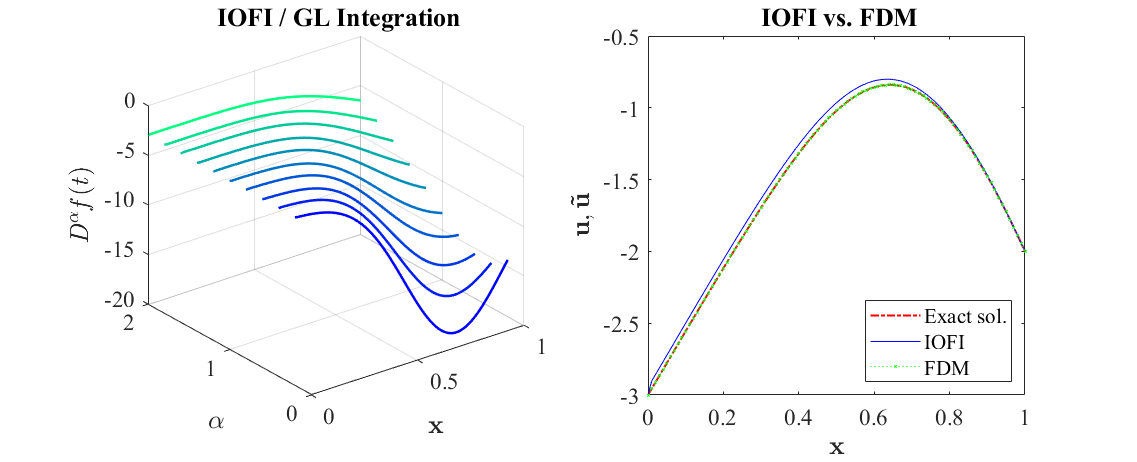}
	\caption{To the left: Evolution of IFOI curves through $x$ and $\alpha$. To the right: Comparison between FDM, IFOI, and the exact solution.}
	\label{fig:caso1}
\end{figure}

The partition on the $\alpha$ domain is regular and rather coarse, with $\Delta \alpha = 0.2$. Partition on $x$ domain is also regular, with $\Delta x = 1.0\edot{-2}$. With these parameters, the method converges successfully by taking $\Delta t_{IOFI} = 1.1\edot{-2}$ s, which is lower than the corresponding value of FDM, $\Delta t_{FDM} = 1.8\edot{-2}$ s. Nonetheless, the error is significantly greater, $e_{IOFI} = 6.2\edot{-2}$ while $e_{FDM}=1.2\edot{-4}$. That is, there are two orders of magnitude of difference. The higher error may be a consequence of the chopping-off of the algorithm used to estimate Eq. \ref{DGL}, as higher-order terms of the series are neglected. 

Figure \ref{fig:caso1} shows results for the first case. To the left, the transformation of the curves as they approach the solution $u_1$, in the sense of Eq. \ref{CLshooting}. The color gradient is related to the growth of $\alpha$. The right-hand side of ODE \ref{caso1} is highlighted in blue, while the solution is shown in green. Notice that the greatest rate of change concerning integration order takes place at the lowest values of $\alpha$. This situation suggests the use of non-regular partitions in this domain to improve error behavior.

Furthermore, Fig. \ref{fig:caso1} shows the deviations of FDM and IFOI approximations with respect to the exact solution to the right. As it is evident, the FDM provides a more accurate answer than IFOI in this case.

\subsection{Case 2}

For this case, the original ODE of case 1 is the same, but mixed boundary conditions are involved, as shown in Eq. \ref{caso2}. Consequently, a variant of the shooting method for Robin's conditions is required. Also, the rectangular product formula replaces the iterative GL-based algorithm in order to estimate the differintegral operator \ref{Dfrac}.

\begin{figure}[ht]
	\centering
	\includegraphics[width=1.0\textwidth]{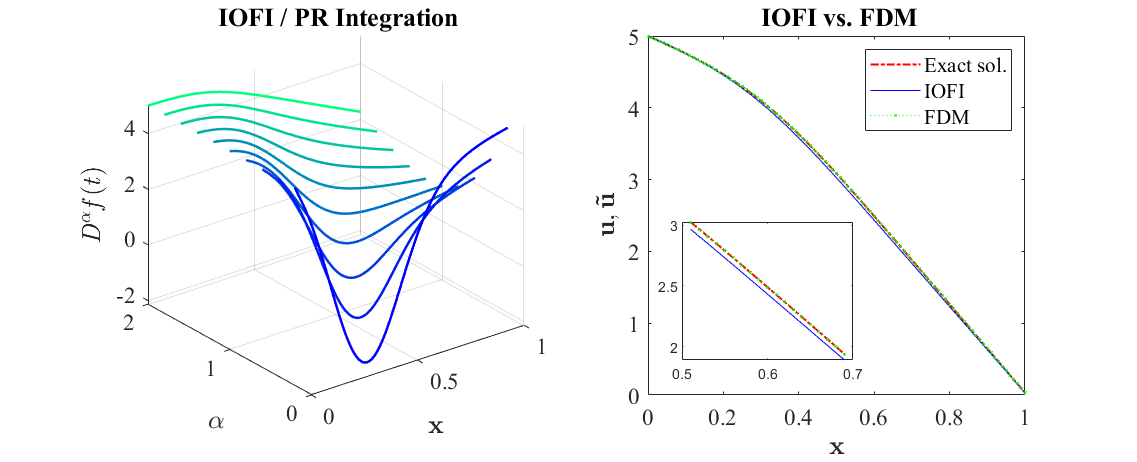}
	\caption{To the left: Evolution of IFOI curves through $x$ and $\alpha$. To the right: Comparison between FDM, IFOI, and the exact solution.}
	\label{fig:caso2}
\end{figure}

Constants are set to $a_2 = 5$, $b_2=200$, and $c_2 = 0.1$ for this case.

Partitions are held identical to those in the former case. As the shooting method converts the original problem into two IVP, the only relevant difference between cases 1 and 2 is the way the fractional operator is approximated.

Analysis shows that although both methods converge, the FDM is faster and more exact than the IFOI method. Specifically, FDM has an execution time of $\Delta t_{FDM} = 1.9\edot{-2}$ s, against $\Delta t_{IOFI} = 2.9\edot{-2}$ s of IFOI. The difference is greater regarding approximation errors, showing a value of $e_{FDM}= 7.3\edot{-5}$ for FDM, while $e_{IOFI} = 5.9\edot{-2}$ for IFOI. Notice that the IFOI error maintains the same order of magnitude as in case 1. As a result, this estimation method for the operator does not seem to report any advantage.

Fig. \ref{fig:caso2} depicts the curve transformations through the $\alpha$ domain. On the left., attenuation of gradients is maximum at lower values of $\alpha$, just as mentioned before. To the right, one can verify the relatively larger error of the IFOI method.

In both cases above, variations on alpha-partitions did not achieve any significant improvement. These alterations are discarded here to keep the text concise.

\subsection{Case 3}

From the previous cases, it becomes apparent that a higher-order scheme is needed to minimize errors. Also, the IFOI method could benefit from an alternative approximation of the fractional operator. The objective is to speed it up, thus reaching competitivity against the FDM. In this section, I have explored the ABMD method, which belongs to the predictor-corrector family techniques.

The ODE structure of previous cases might have prevented an iterative approximation of the integration order from being effective. Even though the former problem has been convenient to assess some estimation schemes, a new structure may enrich the exploration. In this regard, a case such as that of Eq. \ref{caso3} entails a broader range of frequency components. This characteristic may shed new light on the IFOI method.

As previously observed, the error of the IFOI method seems to be higher at the initial iterations, in the sense of Eq. \ref{iterIOFI}. Consequently, in this case, a non-regular partition over $\alpha$ is preferred, so that the integration order follows a quadratic increment.

A fine enough partition is required to avoid \textit{aliasing} phenomena in the numerical method, as the Nyquist theorem states. Nonetheless, it is still worthwhile to examine different outcomes from diverse sizes of $\Delta x$.  A method of insight gain, resulting from the responses, may overcome the downside of sub-sampling.

\begin{table}[htb]
	\centering
	\begin{tabular}{|c|c|c|c|c|c|c|}
		\hline   & \multicolumn{3}{|c|}{Error, $e$}  & \multicolumn{3}{|c|}{Time, $\Delta t$ [s]}   \\
		\hline  & N=40 & N=80 & N=200  & N=40 & N=80 & N=200 \\
		\hline
		\hline
		FDM & 4.8\edot{-4} & 5.9\edot{-5} & 8.6\edot{-6} & 1.8\edot{-2} & 1.8\edot{-2} & 1.8\edot{-2}  \\
		\hline
		IOFI & 6.1\edot{-5} & 5.7\edot{-5} & 3.5\edot{-5} & 3.7\edot{-1} & 6.6\edot{-1} & 1.6\edot{0}\\
		\hline
	\end{tabular}
	\caption{Performance parameters of the FDM and IFOI results.}
	\label{tab:err1}
\end{table}

\begin{figure}[ht]
	\centering
	\includegraphics[width=1.0\textwidth]{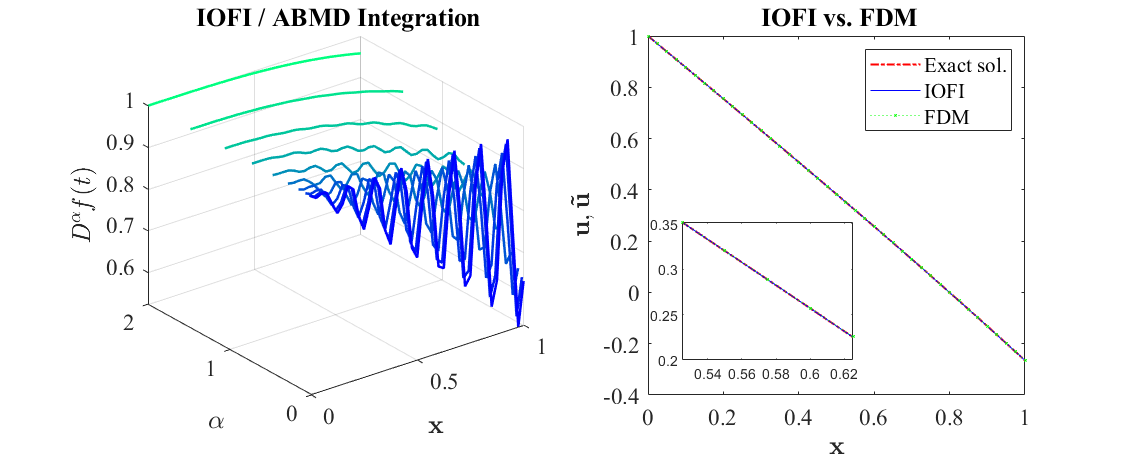}
	\caption{To the left: Evolution of IFOI curves through $x$ and $\alpha$. To the right: Comparison between FDM, IFOI, and the exact solution.}
	\label{fig:caso3}
\end{figure}

Table \ref{tab:err1} presents absolute errors and execution times resulting from each method. Three partition densities are presented, for $N=40$, $N=80$, and $N=200$, with $N$ the number of points over the interval. Out of these values, only the higher $N$ satisfies the Nyquist condition. In this case, both numerical methods have similar errors. Further, the IFOI technique is robust, as its error dependency on $N$ is relatively weaker. On the opposite, the FDM changes in orders of magnitude with each value of $N$. However, from a computational standpoint, execution times reveal that FDM is still substantially more scalable over finer meshes.

Fig. \ref{fig:caso3} illustrates the corresponding graphs to the case $N=40$. To the left, the evolution of the non-homogeneous solution curves, that is $u_1$ from Eq. \ref{CLshooting}. The representation shows the quadratic spacing of the discretization on $\alpha$. This kind of mesh is more effective regarding error treatment than a regular mesh. To the right, the graph shows a comparison between each method's final solutions. Results from both techniques have good quality.

\subsection{Case 4}

At this point, it became interesting to analyze a more complex problem. As mentioned before, case 4 corresponds to a Sturm-Liouville problem type, subject to Dirichlet boundary conditions. The reference method is a variant of the FDM for nonlinear problems. Its implementation follows from the work of \cite{Burden}.

To solve this case a regular partition was set both on $x$ and $\alpha$, with linear increments of  $\Delta x = 0.02$ and $\Delta \alpha = 0.2$. Problem parameters are $a_4 = 3$, and $b_4 = -2$.

\begin{figure}[ht]
	\centering
	\includegraphics[width=1.0\textwidth]{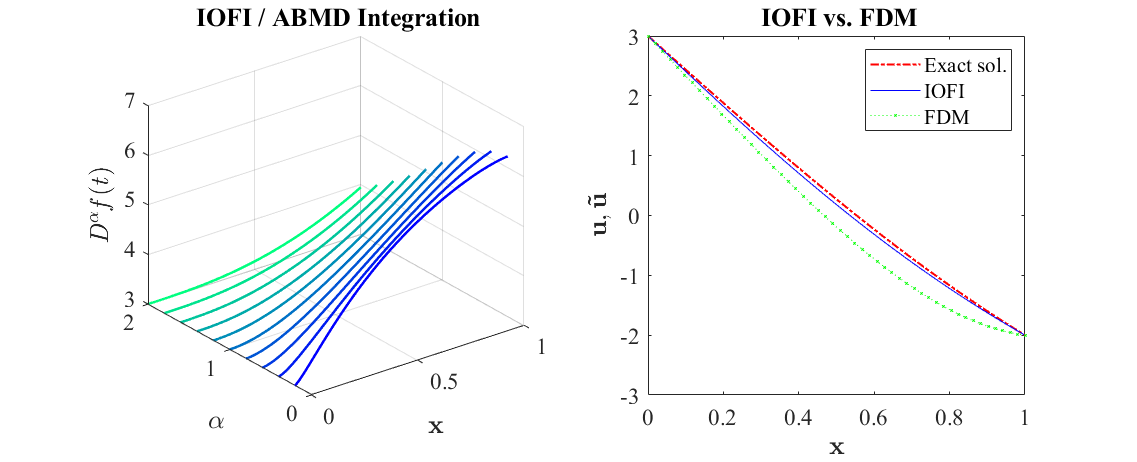}
	\caption{To the left: Evolution of IFOI curves through $x$ and $\alpha$. To the right: Comparison between FDM, IFOI, and the exact solution.}
	\label{fig:caso4}
\end{figure}

The IFOI method turns out to overcome the FDM for this case, as can be seen in Fig. \ref{fig:caso4}. Even though computing times are very similar, with $\Delta t_{IOFI} = 3.5\edot{-1}$ s, and $\Delta t_{FDM} = 3.4\edot{-1}$ s; there is a noticeable difference in each method's error. While the IFOI renders $e_{IOFI}=9.0\edot{-2}$, the reference gives $e_{FDM} = 4.9\edot{-1}$. None of the methods report a significant error sensitivity regarding the mesh diameter.

The graph on the left in Fig. \ref{fig:caso4} shows a progressive evolution to the final solution, distributed regularly for each value of $\alpha$. This fact alone would justify the use of a regular partition for this case. Actually, quadratic spacing, like the one used in the third case, makes the IFOI diverge. The graph on the right allows appreciate the better approximation that IOFI reaches regarding the exact solution compared to the FDM.

\section{CONCLUSIONS}

This article explores an application of fractional calculus to ordinary differential equations by using the semigroup property \ref{semigrupo}. The principal concept is to iterate over the integration-order to approach progressively towards the solution. This concept is named IFOI along with this work.

Precedent analyses enabled a preliminary understanding of the method, together with its strengths and weaknesses. Specifically, the first two cases suggest that the IFOI method is sensitive to chop-off errors. These errors are associated with the numerical approximation method. Ultimately, the existence of error is unavoidable due to the non-local nature of fractional operators. However, higher-order approximation formulae allow an error reduction, closing the gap regarding the reference's error. The ABMD scheme was used for this purpose in cases 3 and 4.

Error sensitivity to truncation can follow from the iterative nature of the IFOI method. Each iteration over the approximation formula amplifies the error. This amplification is repeated for every value of $\alpha$ until it reaches the final problem order. Such a feature becomes a distinctive difference from the classical ODE solution methods, where the approximation formula is required only once to obtain the final solution. Nonetheless, discretization on the $\alpha$ domain shows some advantages. For instance, case 3 reveals certain error insensitivity to mesh diameter compared to the FDM. On the other hand, case 4 has established that the IFOI method can be significantly more effective than FDM for equal mesh diameter.

The evidence from studied cases can not discard the initial hypothesis, which states that it is possible to use coarser partitions in the function domain through IFOI methods while keeping the error equal to the reference method. Even though the computational cost remains a drawback, it seems to be solvable by balancing the order partition with the evolution of the problem dependency concerning integration order.

Further research is needed to fully characterize the IOFI method, both from theoretical and practical viewpoints. Due to its nature, this technique may become adequate to solve strongly nonlinear, stiff, or ill-conditioned problems. Future work will focus on these issues. Moreover, research into solving partial differential equations may be promising. The additional domains of integration orders would enable a substantial optimization of the method due to the higher degrees of freedom.

\bibliographystyle{unsrt}


\end{document}